\documentclass[english,a4paper,12pt,dvipdfm]{amsart}
\usepackage[all,web]{xy}
\usepackage{graphicx}
\usepackage{amsmath} 
\usepackage{amsfonts}
\usepackage{amssymb}
\usepackage{amsthm}
\usepackage{setspace}

 \theoremstyle{plain}    
 \newtheorem{thm}{Theorem}[section]
 \numberwithin{equation}{section} 
 \numberwithin{figure}{section} 
 \theoremstyle{remark}   
  
 \theoremstyle{plain}    
 \theoremstyle{plain}    
  
 \theoremstyle{remark} 
 
 \theoremstyle{remark}
 
 \theoremstyle{definition}
  
 \theoremstyle{plain}
 \theoremstyle{plain}    
 \theoremstyle{plain}    
 \theoremstyle{definition}
 
 \theoremstyle{definition}
  
 \theoremstyle{plain}    
  
 \theoremstyle{plain}    
 \theoremstyle{remark}    
  
 \theoremstyle{remark}    
  
 \theoremstyle{definition}
  
 \theoremstyle{plain}
 \theoremstyle{remark}
 \newtheorem{rem}[thm]{Remark}
 \theoremstyle{remark}

 1

\def\CC{\mathbb{C}}

\def\PP{\mathbb{P}}
\def\QQ{\mathbb{Q}}
\def\RR{\mathbb{R}}

\def\ZZ{\mathbb{Z}}

\newcommand{\mtwo}[4]{\left(
        \begin{matrix}#1&#2\\#3&#4
        \end{matrix}\right)}

\def\inner<#1,#2>{{\left\langle{{#1},{#2}}\right\rangle}}
\def\sl2of#1{\textrm{SL}_2(#1)}

\def\inner<#1,#2>{{\left\langle{{#1},{#2}}\right\rangle}}

\usepackage{babel}

\begin{document}

\title{A note on the rational points of $X_0^+(N)$}
\author{Carlos Casta\~no-Bernard}
\email{ccastanobernard@gmail.com}

\begin{abstract}
Let $C$ be the image of a canonical embedding $\phi$ of
the Atkin-Lehner quotient $X_0^+(N)$
associated to the Fricke involution $w_N$.
In this note we exhibit some relations among the 
rational points of $C$.
For each $g=3$ (resp. the first $g=4$) curve $C$
we found that there are one or more
lines (resp. planes) in ${\bf P}^{g-1}$ whose intersection with $C$ consists
entirely of rational Heegner points or the cusp point,
where $N$ is prime.
We also discuss an explanation of the first non-hyperelliptic exceptional
rational point.
\end{abstract}


\maketitle

\pagenumbering{roman}
\setcounter{page}{0}



\tableofcontents


\pagenumbering{arabic}
\pagestyle{headings}

%
%
%
\section{Introduction}
Fix an integer $N>1$ and
let $X_0(N)$ be the moduli space of
(ordered) pairs $(E,E^\prime)$ of
generalised elliptic curves $E$ and $E^\prime$ linked by
a cyclic isogeny $\varphi\colon E\longrightarrow E^\prime$
of degree $N$.
Consider the Atkin-Lehner quotient curve $X_0^+(N)$
defined by
the involution $w_N$ of $X_0(N)$ induced by mapping
an isogeny $\varphi\colon E\longrightarrow E^\prime$
to its dual $\hat{\varphi}\colon E^\prime\longrightarrow E$.
The quotient curve $X_0^+(N)$ has been studied by
Galbraith~\cite{galbraith:rational},
Mazur~\cite{mazur:eisen},
and Momose~\cite{momose:xplus},
among others. 
Galbraith~\cite{galbraith:rational} studied the
rational points of a canonical image $C\subset\PP^{g_N^+ -1}$ of $X_0^+(N)$,
where $g_N^+$ is the genus of $X_0^+(N)$.
In particular he exhibits explicit formul\ae\ for $C$ for
prime conductors $N$ for which the curve has genus $g\leq 5$.
In each case he
locates the cusp and rational CM points and,
moreover,
for $N=137$ he exhibits a rational point which is
neither a cusp point nor a CM point. 
In this note we exhibit an explicit set of
hyperplanes $\{H_1,\dots H_s\}$ in $\PP^{g_N^+ -1}$
such that the intersection of each $H_i$ with $C$ (over the complex numbers)
consist entirely of rational points of $C$,
for each prime level $N$ such that $g_N^+=3$,
i.e. $N=97$, $109$, $113$, $127$, $139$, $149$, $151$, $179$, and $239$,
and the first prime level $N$ such that $g_N^+=4$,
i.e. $N=137$.
For the latter case we found a further plane
defined by $3$ different CM points that contains the exceptional point.

The material is organised as follows.
Section~\ref{sec:gal} introduces some basic results
that we used to compute an equation for $C$.
The collinearity relations are discussed in
Section~\ref{sub:genus3},
while the coplanarity relations are discussed in Section~\ref{sub:genus4}.

\bigskip
\textsc{Acknowledgments}
I would like to thank the referee of Experimental Mathematics
for his/her insightful comments,
which led to some significant improvements on a previous version of this paper.
I would also like to thank Barry Mazur for his interesting remarks.
I would like to heartily thank Jenny and Kenneth Cowan
for their hospitality while preparing this paper.

%
\section{Preliminaries}\label{sec:gal}
Let $X$ be an algebraic curve
defined over a field $k$ and let $\Omega^1_X$
be the $k$-vector space of its holomorphic differentials.
Also let $\{\omega_1,\dots\omega_g\}$ be a basis of $\Omega^1_X$.
The integer $g$ is called the \textit{genus} of $X$.
The \textit{canonical map} $\phi$ of $X$
in projective space $\PP^{g-1}$ is the morphism
\begin{equation*}
\xymatrix{
X \ar [r]^{\phi}& \PP^{g-1}\\
P \ar @{|->}[r]& (\omega_1(P):\dots:\omega_g(P))\\
}
\end{equation*}
It is well-known that the canonical map $\phi$ is an embedding,
if the genus $g$ exceeds $2$ and $X$ is not hyperelliptic.
Now fix an integer $N>1$ and recall that
\begin{equation*}
X_0(N)=\Gamma_0(N)\backslash\mathcal{H}^*,
\end{equation*}
where
\begin{equation*}
\Gamma_0(N)=
\left\{\mu=\mtwo{\alpha}{\beta}{\gamma}{\delta}\in\textrm{SL}_2(\ZZ)
\,\colon\,
\gamma\equiv 0\pmod{N}\right\}.
\end{equation*}
The Atkin-Lehner quotient curve $X_0^+(N)$ may be defined as
the curve $X_0(N)$ modulo de action of the Fricke involution $w_N$,
which is defined by
\begin{equation*}
\tau\mapsto-\frac{1}{N\tau}.
\end{equation*}
Every holomorphic differential $\omega$ on $X_0(N)$
may be expressed uniquely as $\omega=fd\tau$,
where $f(\tau)=\sum_{n=1}^\infty a(n)q^n$ (where $q=e^{2\pi\tau}$)
is a modular form for $\Gamma_0(N)$ and of weight $2$  ,
and conversely every modular form for $\Gamma_0(N)$ and of weight $2$
gives rise to a holomorphic differential $\omega=fd\tau$.
In particular the canonical maps for $X_0^+(N)$ are of the form
\begin{equation*}
\xymatrix{
X_0^+(N) \ar [r]^{\phi}& \PP^{g-1}\\
P \ar @{|->}[r]& (f_1(\tau):\dots:f_{g_N^+}(\tau)),\\
}
\end{equation*}
where $\{f_1,\dots f_{g_N^+}\}$ constitute a basis of
the $+1$-eigenspace $S_2^+(\Gamma_0(N))$
of the $\CC$-vector space of
modular forms $S_2(\Gamma_0(N))$ with respect to the action of
the Fricke involution $w_N$.
It is a classical fact that the vector space $S_2(\Gamma_0(N))$
has a basis consisting of modular forms
with rational integer coefficients only.
The same is true for the eigenspace $S_2^+(\Gamma_0(N))$.
A set of equations with integral coefficients $S$
for the curve $C$ may be computed by
finding combinations of powers of the $q$-expansions
which yield identically zero series.
Using the method described in Galbraith~\cite[p. 19]{galbraith:thesis},
with the help of William Stein's
computer package {\sc Hecke}~\cite{stein:hecke}
and \textsc{Pari}~\cite{pari:gp}
it is easy to compute these equations $S$ for small genera,
such as $g_N^+=3$ and $g_N^+=4$.
(See also Elkies~\cite{elkies:ell} for related methods and results on the
computation of equations for modular curves.)
In these two cases Galbraith~\cite{galbraith:thesis} has shown that
the curves $X_0^+(N)$ are non-hyperelliptic.
So in fact $C$ is a complete intersection.
(See Examples IV.5.2.1 and IV.3.3.2 of Hartshorne~\cite{hartshorne:alg}.)
Once we have a set of such equations $S$ for $C$ in practice it only takes
only a brute-force search to obtain as many rational points on $C$
as predicted by the theory of Complex Multiplication
(see Gross~\cite{gross:heegner}),
the cusp rational point,
or some other rational points that come as
a fixed point of certain hyperelliptic involution.
Sometimes there are more rational points than these ``obvious'' ones;
Galbraith~\cite{galbraith:thesis} has shown that there exist
rational points in $C$ that are neither cusps nor
CM rational points for non-hyperelliptic $X_0^+(N)$ of genus at least $4$.
 
From now on let us assume that $N$ is prime.
Then using Proposition 3.1 of Gross~\cite[p. 347]{gross:primelevel}
and the Riemann-Hurwitz formula we may see that
\begin{displaymath}
g_N^+ = \frac{1}{2}\left(g_N + 1 - H(N)\right),
\end{displaymath}
where
\begin{displaymath}
H(N)=\left\{
\begin{array}{ll}
\frac{1}{2}h(-4N),& \textrm{if $N\equiv 1 \pmod{4}$}\\
\frac{1}{2}(h(-N) + h(-4N)),& \textrm{otherwise,}\\
\end{array}
\right.
\end{displaymath}
with $h(D)$ the class number of
the imaginary quadratic order of discriminant $D$,
and $g_N$ is genus of $X_0(N)$,
which is given by $g_N = \lfloor\frac{N+1}{12}\rfloor$,
unless $N=12q+1$ when $g=q-1$.
In particular,
by using explicit upper bounds on the class number $h(D)$
it may be found that the prime conductors $N$ such that $X_0^+(N)$ has genus $3$
are indeed $N=97$, $109$, $113$, $127$, $139$, $149$, $151$, $179$, and $239$.
and all these curves $X_0^+(N)$ are non-hyperelliptic.
Similarly it may be found that
there are exactly $5$ prime numbers $N$
such that $X_0^+(N)$ has genus $g=4$,
namely $N = 137$, $173$, $199$, $251$, and $311$.

%
\section{Genus three: collinearity relations}\label{sub:genus3}
Let us assume $N$ is one of the $9$ prime conductors $N$
such that $X_0^+(N)$ has genus $g_N^+=3$.
So the canonical map $\phi$ associated to
a basis of modular forms with
integral Fourier coefficients $\{g_1,g_2,g_3\}$ is
an embedding $\phi$ defined over $\QQ$
of $C$ into the projective plane $\PP^2$
such that its image $C$ is of degree $4=2(g_N^+-1)$.
So we may obtain a projective equation $F(X, Y, Z)=0$ for the plane curve $C$
by computing a (non-trivial) linear relation among the elements in the set
\begin{equation*}
\{g_1^ag_2^bg_3^c\in S_8(\Gamma_0(N))\colon a,b,c\in\ZZ_{\geq 0}, a+b+c=4\}.
\end{equation*}
Note that a line $L$ in the projective plane $\PP^2$ will intersect
the curve $C$ in $4$ points,
if we take into account intersection multiplicities.
So heuristically,
the line defined by two rational points $P_1$ and $P_2$ of $C$
(which is the tangent line of $C$ at say $P_1$, in case $P_1=P_2$)
is not expected to intersect $C$ in further rational points;
the remaining points of intersection are expected to be in general
two quadratic irrationals (one conjugate to the other).
However,
as we shall see it is possible to
exhibit for each of the levels under consideration
a non-empty set of lines $\{L_1,\dots L_s\}$
such that the intersection of each $L_i$ with $C$
(over the complex numbers) consists entirely of rational points.
We describe these lines with the help of some diagrams below.
Each diagram contains an equation for $C$ and a table of rational points.
Since all these rational points are CM points, or the cusp point,
we label these rational points according to the discriminant $D$
of the CM point.
By the work of Galbraith~\cite{galbraith:thesis}
we know that in this case the relevant values of $D$ are
\begin{displaymath}
D = -3,-4,-7,-8,-11,-12,-16,
\end{displaymath}
\begin{displaymath}
-19,-27,-28,-43,-67,-163.
\end{displaymath}
To ease notation we regard the cusp $i\infty$ labeled by $D=0$.
The discriminant $D$ of
each rational point of $C$ exhibited was recognized
in most cases by computing a suitable approximation of the image in $C$ of
the Heegner point of $X_0^+(N)$ of discriminant $D$
so that it would be clear which of the rational points of $C$
found previously corresponds to $D$.

\begin{rem}
Sometimes simply too many Fourier coefficients were necessary to
obtain a suitable approximation
for the image in $C$ of a rational Heegner point,
e.g. for $N=97$.
In these cases there were exactly two rational points of $C$,
say $P$ and $Q$,
which needed further work to be labeled.
We were able to resolve these kind of ambiguities by using a variant of
a result of Ogg~\cite{ogg:real} which may be applied to 
describe the real locus $X_0^+(N)(\RR)$ in terms of certain
indefinite binary quadratic forms.
More precisely,
first we were able to locate (exactly) in
a parametrisation $\pi$ of $X_0^+(N)(\RR)$
all the Heegner (i.e. CM) rational points,
and then by ``walking along''
a good polygonal approximation of the real locus of $C$
in the direction determined by $\pi$,
we paired the unmatched discriminants,
say $D_1$ and $D_2$,
with the points $P$ and $Q$.
\end{rem}

The diagrams are self explanatory;
the intersection multiplicity $\partial$ of $L_i$ with $C$
at a specific point $P\in L_i\cap C$ is indicated only when necessary.

\begin{center}
\includegraphics[scale=0.8]{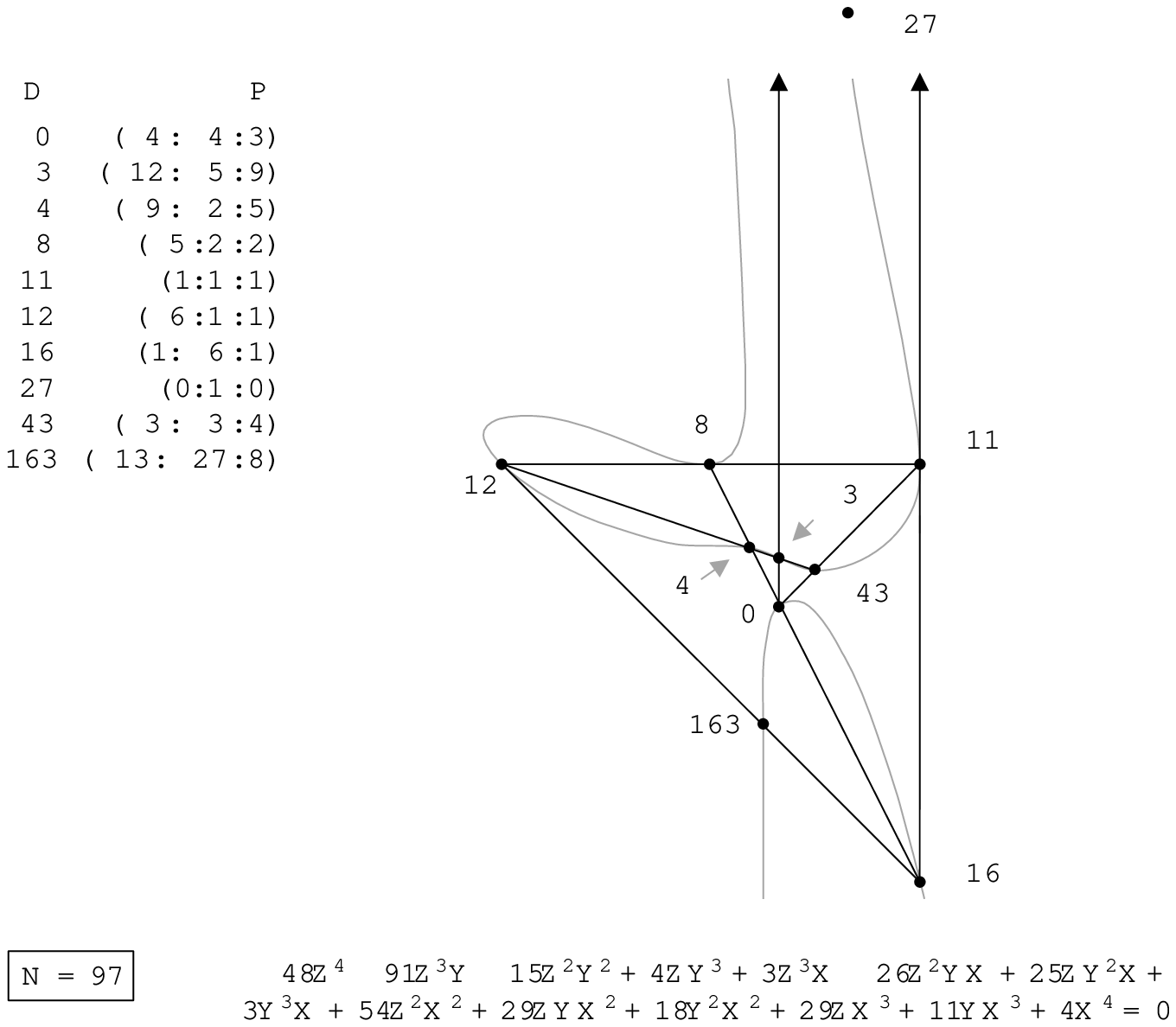}%
\end{center}

\begin{center}
\includegraphics[scale=0.8]{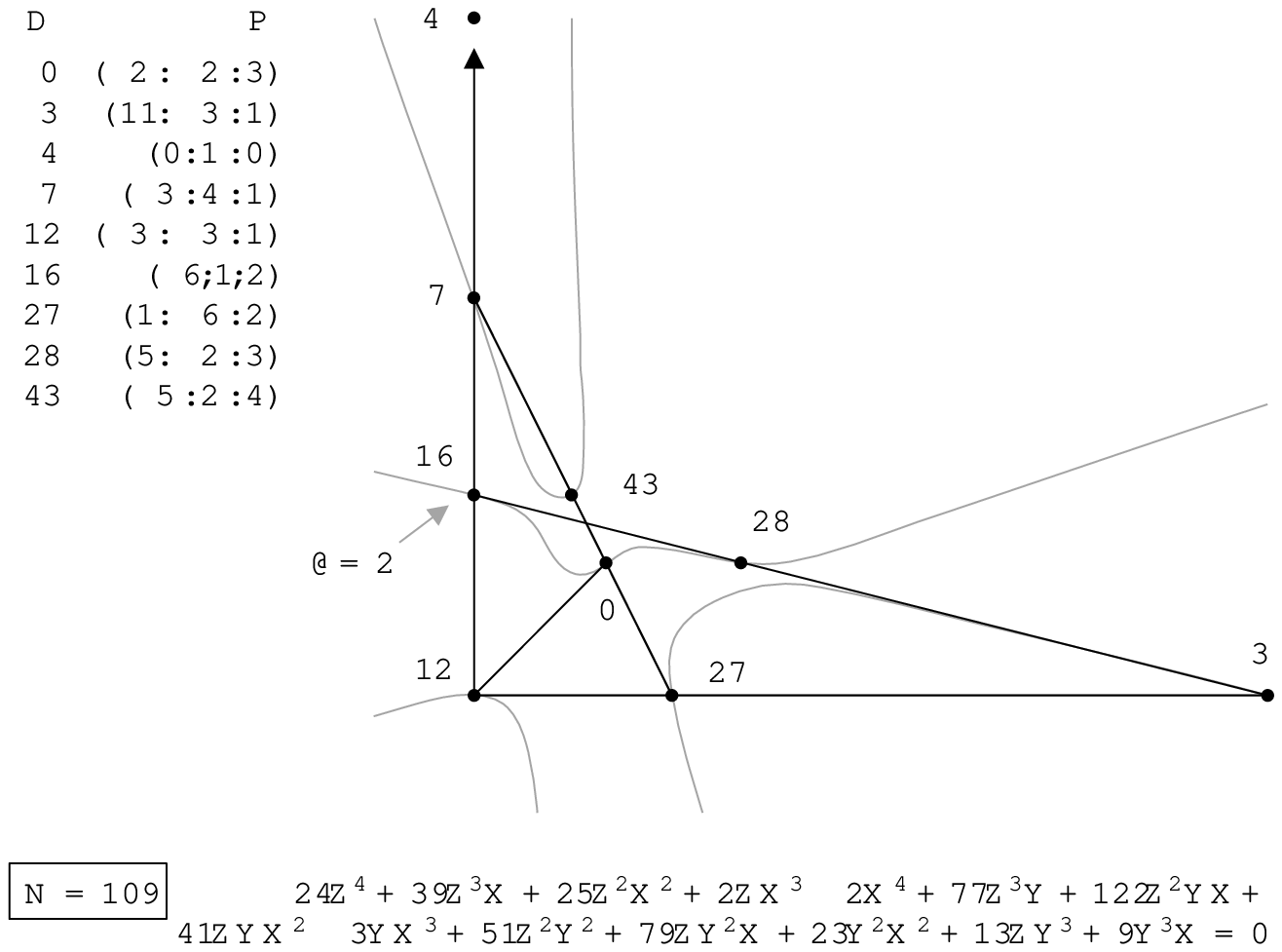}%
\end{center}

\begin{center}
\includegraphics[scale=0.8]{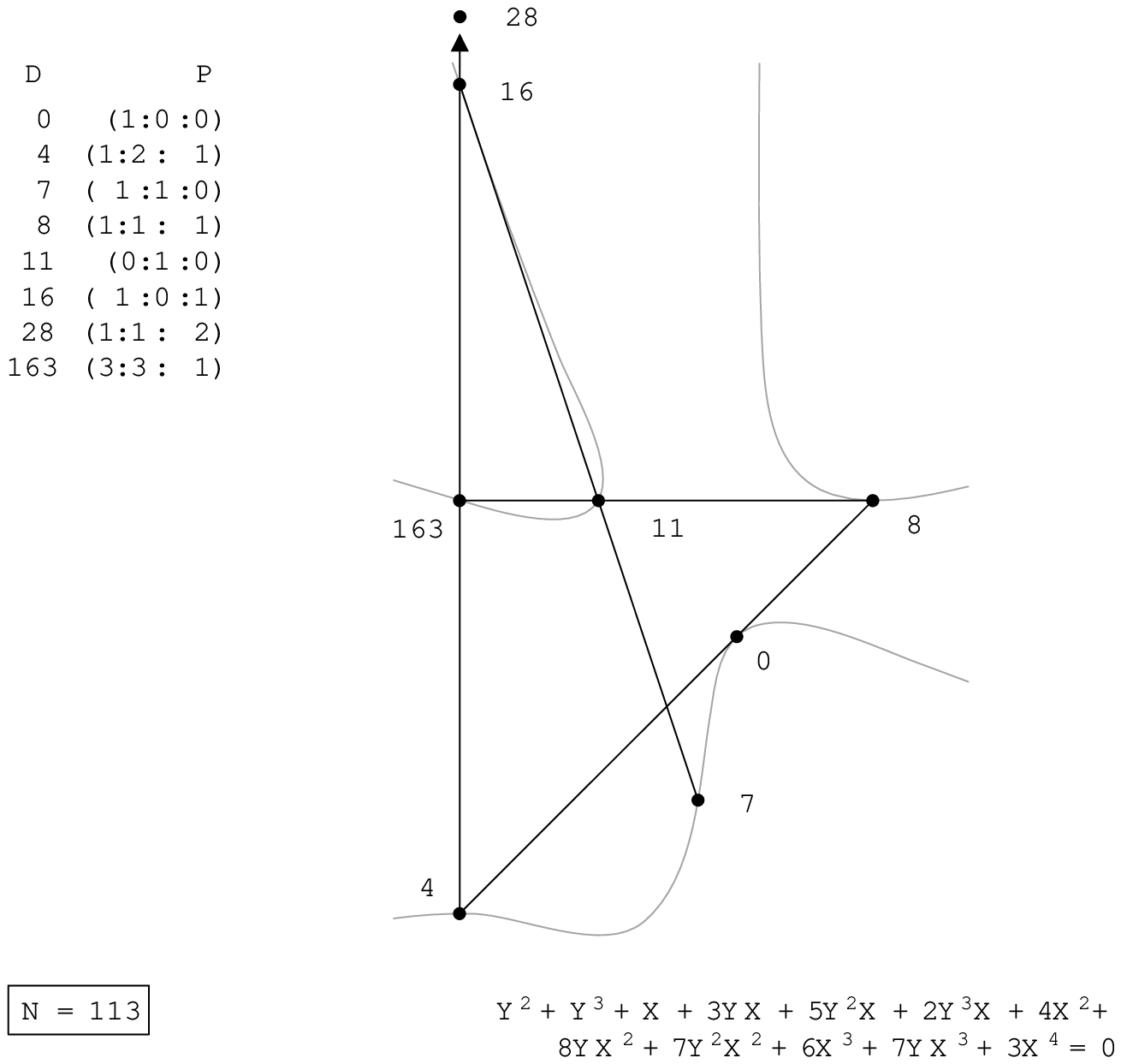}%
\end{center}

\begin{center}
\includegraphics[scale=0.8]{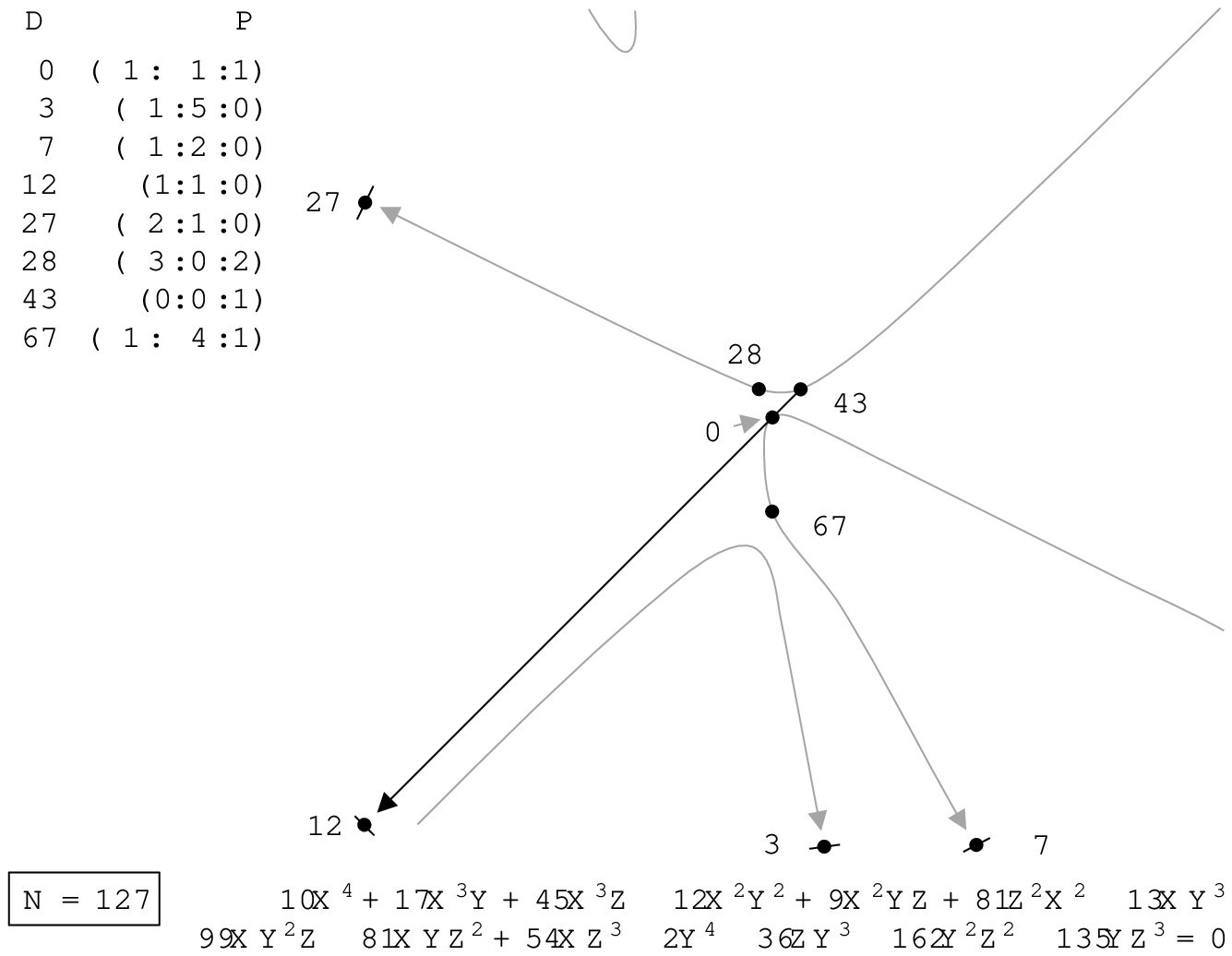}%
\end{center}

\begin{center}
\includegraphics[scale=0.8]{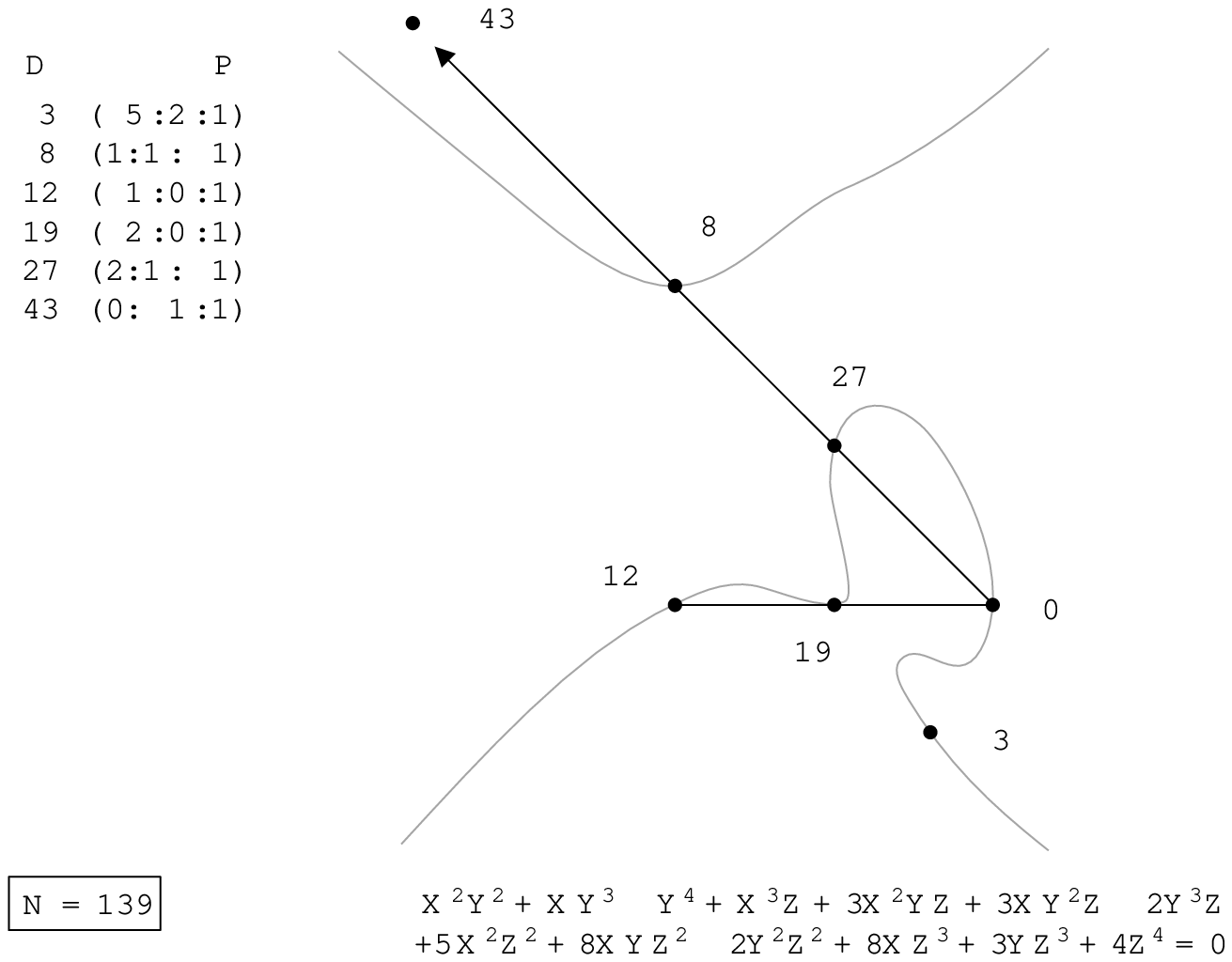}%
\end{center}


\begin{center}
\includegraphics[scale=0.8]{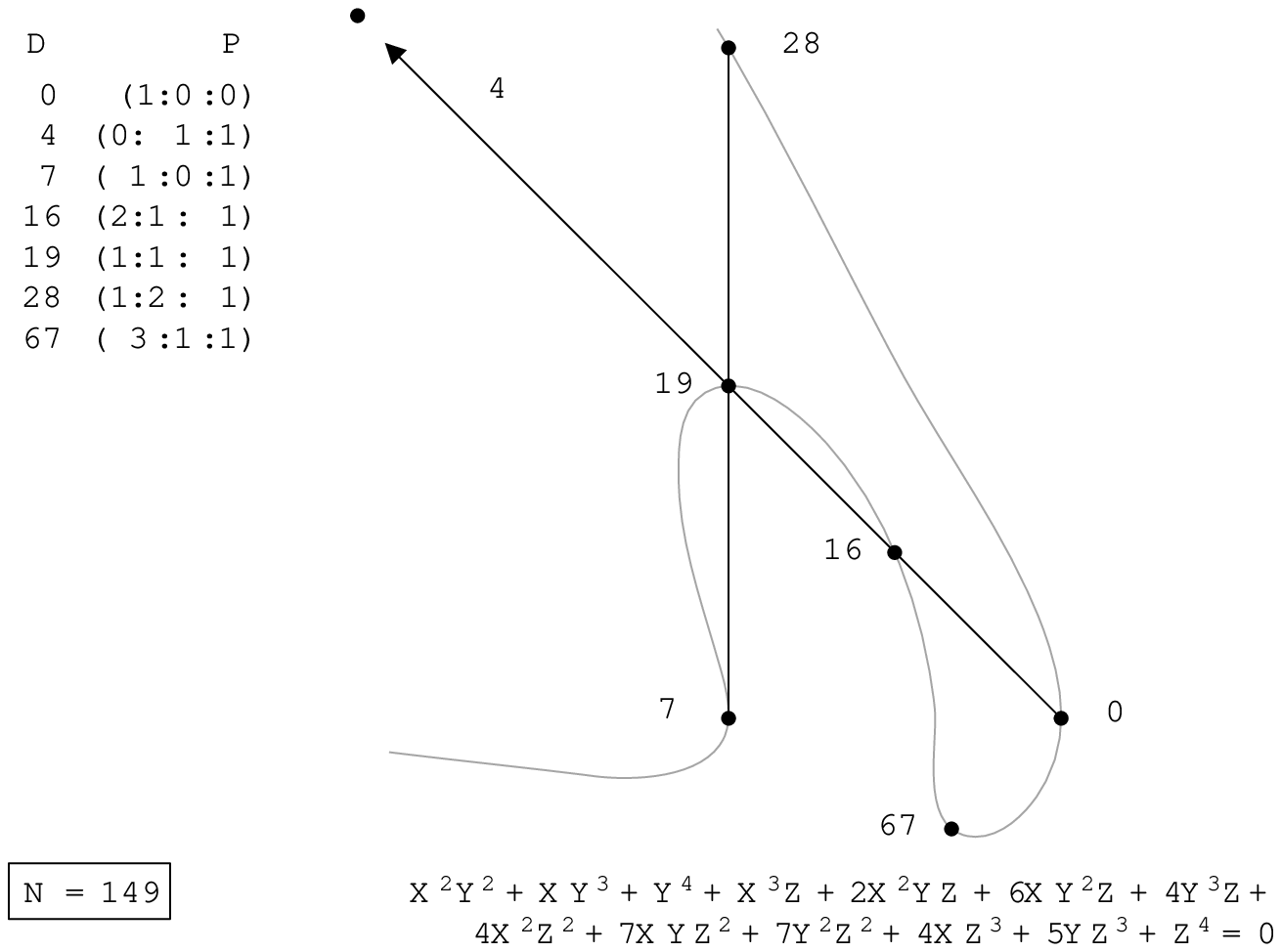}%
\end{center}

\begin{center}
\includegraphics[scale=0.8]{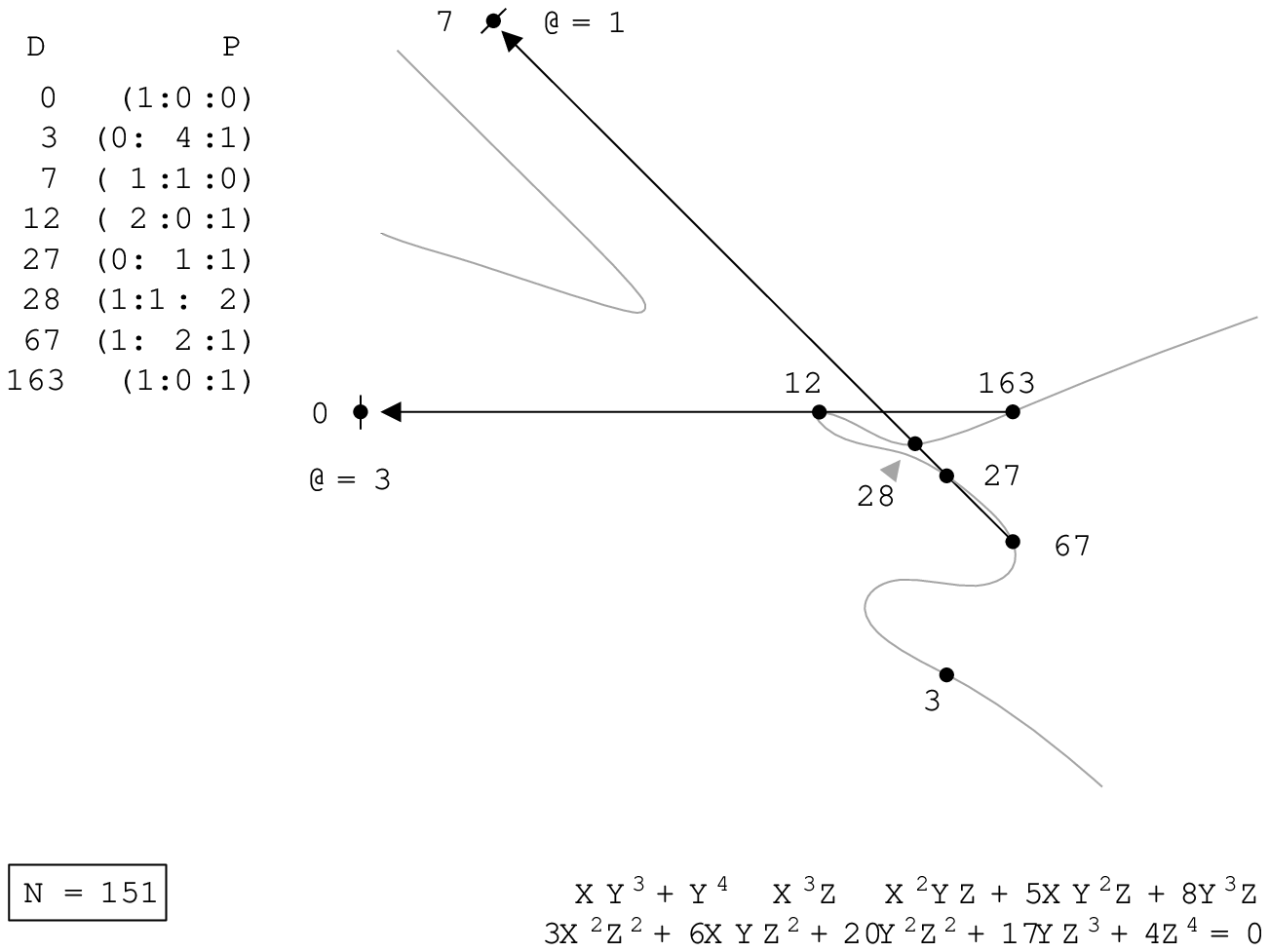}%
\end{center}

\begin{center}
\includegraphics[scale=0.8]{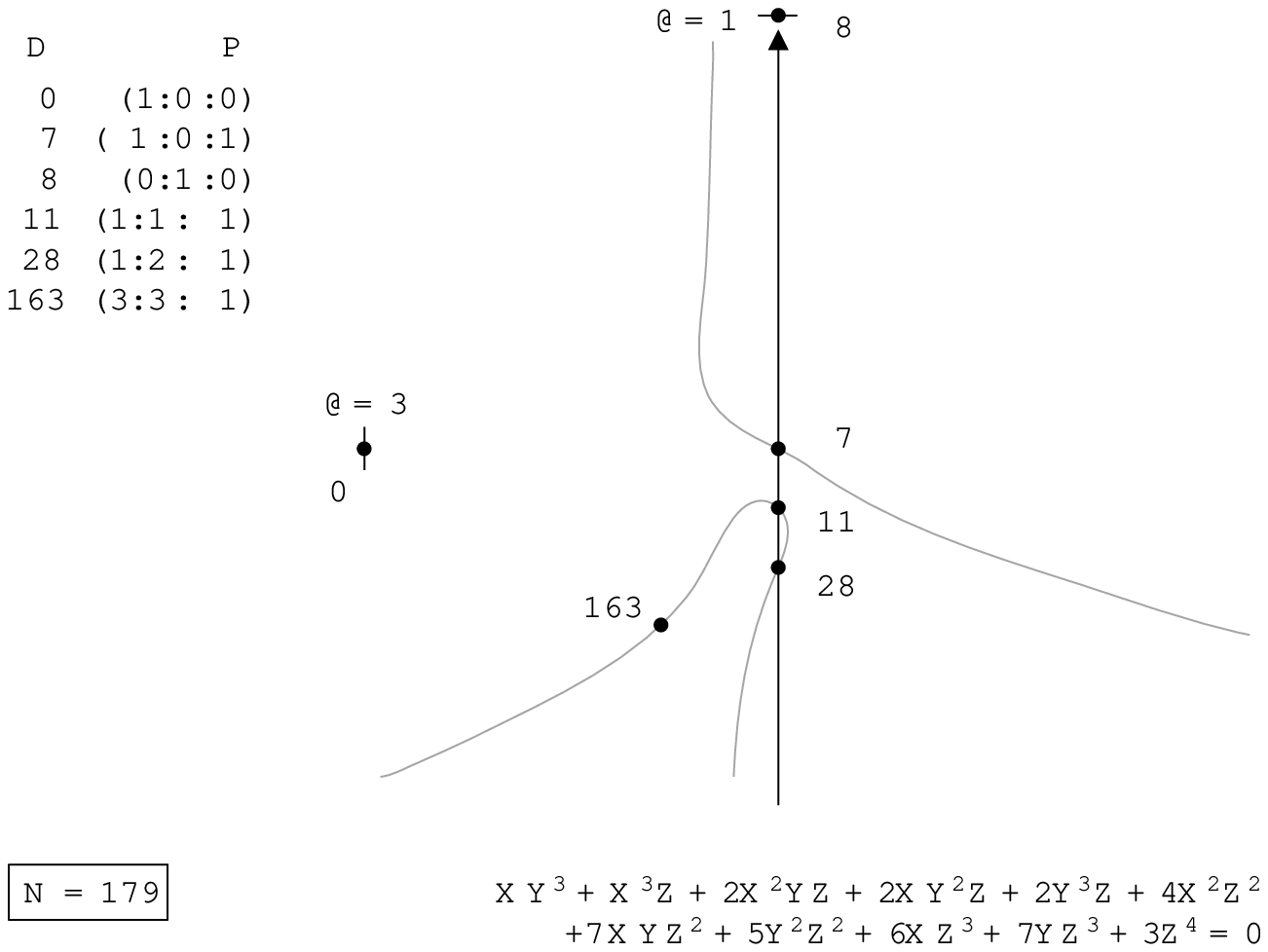}%
\end{center}

\begin{center}
\includegraphics[scale=0.8]{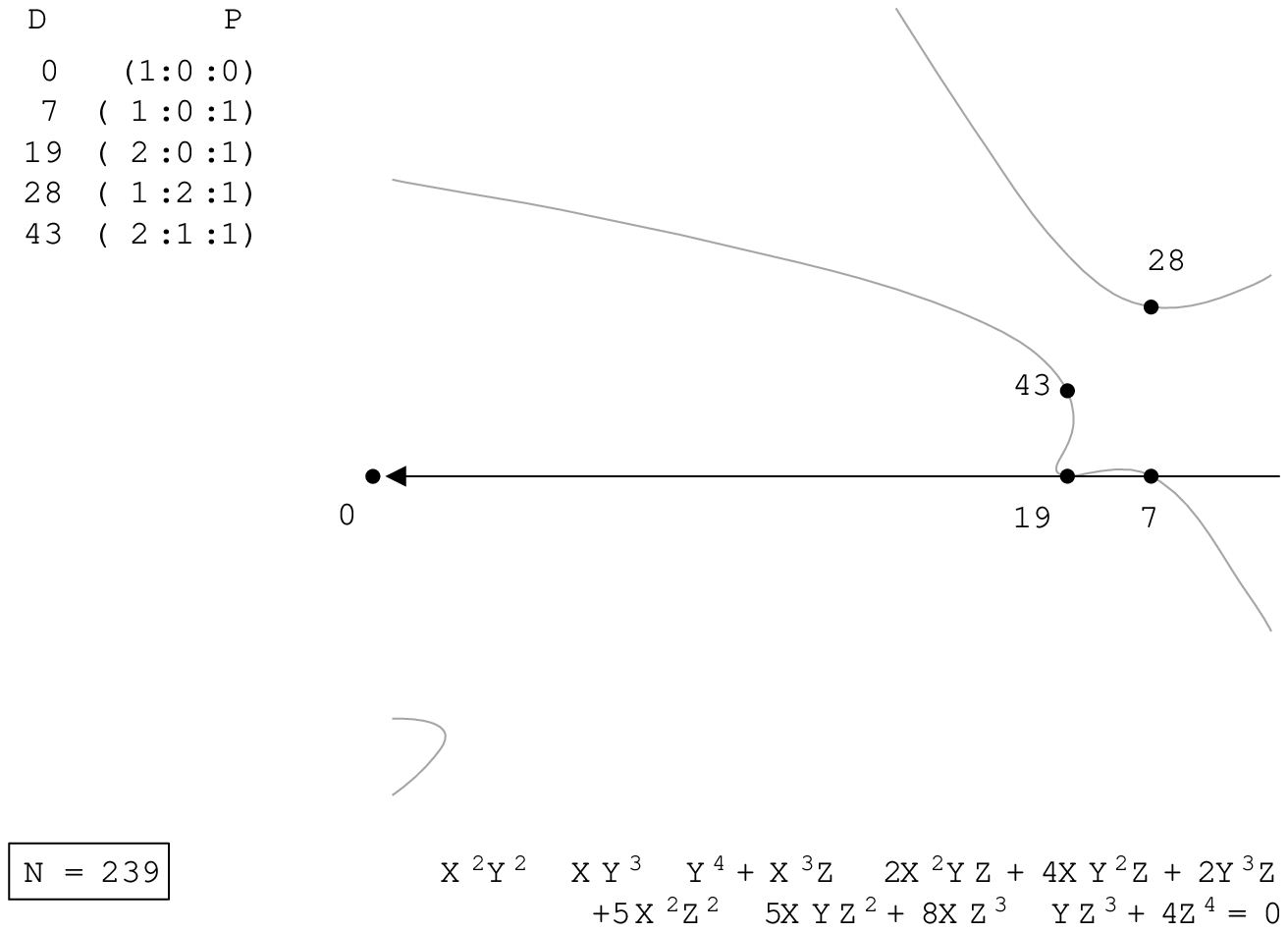}%
\end{center}

%
\section{Genus four: coplanarity relations}\label{sub:genus4}
Now we will study the rational points of
the first genus $4$ Atkin-Lehner quotient curve $X_0^+(N)$.
Again by the work of Galbraith~\cite{galbraith:thesis}
we know that $C$ may be given as the intersection in $\PP^3$ of
\begin{equation*}
XY + WY + 2Y^2 + 2WZ + XZ + 6YZ + 3Z^2 = 0,
\end{equation*}
and
\begin{equation*}
\begin{split}
X^3 + WX^2 + 6X^2Z - 2XY^2 - 5XYZ + XZW + 13XZ^2 + 2Y^3 \\
+ 3WY^2 + W^2Y + 3WYZ - 6YZ^2 + ZW^2 - 4Z^2W + 14Z^3 = 0.
\end{split}
\end{equation*}
(See also Example IV.5.2.2 of Hartshorne~\cite{hartshorne:alg}.)
By a brute-force search we may find in $C$ the rational points
\begin{displaymath}
\begin{array}{rr}
D & P\\
\vspace{-2ex}\\
0&[1:0:0:0]\\
-4&[2:-4:-3:2]\\
-7& [2:-1:-2:1]\\
-8& [-1:1:0:0]\\
-11& [1:1:-1:0]\\
-16& [2:0:-1:0]\\
-19& [1:-2:-1:1]\\
-28& [0:1:2:-1]\\
\end{array}
\end{displaymath}
We may also find in $C$ the non-CM rational point $[19:2:-16:4]$.

The fact that the degree of $C$ is $6$
implies that the planes $\Pi$ in $\PP^3$
will intersect $C$ (over a fixed algebraic closure of $\QQ$) at $6$ points,
if we take into account multiplicities.
So if $\Pi$ is a plane defined by $3$ rational points of $C$,
then the other intersection points of $\Pi$ with $C$
are in general expected to be defined over a cubic extension of $\QQ$.
Hovever,
it turns out that
there are $3$ different planes $\Pi_1$, $\Pi_2$, and $\Pi_3$ in $\PP^3$
such that each of these planes intersects $C$ at exactly $6$ rational points
(with multiplicities):
\begin{itemize}
\item $\Pi_1$: $z = 0$: $$2(0) + 2(-8) + (-11) + (-16)$$
\item $\Pi_2$: $x + 2y + 3z + 1 = 0$ $$(-7) + (-8) + 2(-11) + (-16) + (-19)$$
\item $\Pi_3$: $x + y + 3z = 0$ $$(0) + 2(-7) + (-11) + (-19) + (-28)$$
\end{itemize}
There are some other remarkable properties about the set of rational points
of $C$.
The set $\{-7, -11, -19\}$ is contained in one line, say, $L_1$,
and also the set $\{-8, -11, -16\}$ is contained in one line, say, $L_2$.
These two lines meet at $D = -11$ and,
moreover,
both are contained in plane $\Pi_2$.
Finally,
the $3$ planes $\Pi_1$, $\Pi_2$, $\Pi_3$ meet at $D=-11$..

Another property worth mentioning is that the plane $\Pi_e$ defined by
the points $D=0,-4,-11$
\begin{displaymath}
2x + 2y + 7z = 0
\end{displaymath}
also contains the exceptional point,
which is quite unexpected since this may not be heuristically explained by
the small size of the coefficients of the Jacobi form involved.
(See Gross, Kohnen and Zagier~\cite{gross:gkz}.)
However,
the remaining two points of intersection are not rational;
these are two conjugate points defined over the real quadratic field of
discriminant $\Delta=8$.

\section{Concluding remarks}
It seems to be an open question
whether these relations may have an interesting explanation,
or else if these relations are only
a consequence of an ``accident'' due to
the small size of the genus $g_N^+$ of $X_0^+(N)$,
perhaps related to the fact that the field $K_f$ generated
by the Fourier coefficients of the newform $f$ of level $N$
and ``$-1$'' sign in the functional equation of
its $\Lambda$-function,
is an abelian extension of $\QQ$ for each
prime $N$ with $g_N^+=3$.
It seems worthwhile to
extend the above list of examples to higher levels,
hoping that a more
extensive experimental evidence
may help to grasp the nature of this phenomenon.
This may shed some more light into
the nature of $X_0^+(N)(\QQ)$ for prime levels $N$,
which is ``extremely interesting'',
as expressed in Mazur~\cite{mazur:eisen}.

%
%
%
\bibliographystyle{amsplain}
\bibliography{biblio}
\printindex
\end{document}